\documentclass[amstex,12pt]{article}
\usepackage{amssymb,amsmath,amsthm,graphicx}
\newtheorem{theorem}{Theorem}[section]

\newtheorem{lemma}{Lemma}[section]
\newtheorem{definition}{Definition}[section]
\newtheorem{remark}{Remark}[section]
\newtheorem{example}{Example}[section]

\newcommand{\beq}{\begin{equation}}
\newcommand{\eeq}{\end{equation}}
\newcommand{\beqn}{\begin{eqnarray}}
\newcommand{\eeqn}{\end{eqnarray}}

\allowdisplaybreaks \baselineskip 20pt
\begin{document}

\title{Some remarks on  almost periodic time scales and
almost periodic functions on time scales\thanks{This work is supported by the National Natural
Sciences Foundation of People's Republic of China under Grant
11361072.}}
\author{Yongkun Li\thanks{Email: yklie@ynu.edu.cn.}\\
Department of Mathematics,
Yunnan University\\
Kunming, Yunnan 650091\\
 People's Republic of China}

\date{}
\maketitle \allowdisplaybreaks
\begin{abstract}
 In this  note we communicate some important remarks about the concepts of  almost periodic time scales and
almost periodic functions on time scales that are proposed by  Wang and Agarwal  in their recent papers (Adv. Difference Equ.  (2015) 2015:312;  Adv. Difference Equ.  (2015) 2015:296;  Math. Meth. Appl. Sci. 2015, DOI: 10.1002/mma.3590).
\end{abstract}

\allowdisplaybreaks
\section{Introduction}

 \setcounter{section}{1}
\setcounter{equation}{0}
 \indent

In order to study  almost periodic functions, pseudo almost periodic functions, almost automorphic functions, pseudo almost automorphic functions and so on, on time scales,  the following concept of almost periodic time scales was proposed in \cite{li1,liwang}:
\begin{definition}\label{pe}\cite{li1}  A time scale $\mathbb{T}$ is called an almost periodic time scale if
\begin{eqnarray}\label{epe}
\Pi=\big\{\tau\in\mathbb{R}: t\pm\tau\in\mathbb{T}, \forall t\in{\mathbb{T}}\big\}\neq\{0\}.
\end{eqnarray}
\end{definition}

Based on Definition \ref{pe}, almost periodic functions \cite{li1,liwang}, pseudo almost periodic functions \cite{liw1}, almost
automorphic functions \cite{lm}, almost automorphic functions of order $n$ \cite{mo}, weighted  pseudo almost
automorphic functions \cite{liw2} and weighted piecewise pseudo almost automorphic functions \cite{wa5} on time
scales were defined successfully. 

Since Definition \ref{pe} requires that time scale $\mathbb{T}$ has a global property, that is, there exists at least one $\tau\in \mathbb{R}$ such that $t\pm\tau\in \mathbb{T}$ for all $t\in \mathbb{T}$,
it is very restrictive. This may exclude  many interesting time scales. Therefore, it is a
challenging and important problem in theory and applications to find   new concepts of almost periodic time
scales.  Recently, Wang and Agarwal  in \cite{wa1,wa2,wa4} have made some efforts to introduce some new types of almost periodic time scales and almost periodic functions on time scales. However, unfortunately, there are many flaws and mistakes in \cite{wa1,wa2,wa4}.

Our main purpose of this note is to point out some flaws and mistakes in \cite{wa1,wa2,wa4} and give another correction of the concept of almost periodic functions on time scales in \cite{li1}, which is overcorrected in \cite{wa1}.

\section{ Preliminaries}

  \setcounter{equation}{0}

\indent

In this section, we shall first recall some definitions and state some results which are used in what follows.

Let $\mathbb{T}$ be a nonempty closed subset (time scale) of
$\mathbb{R}$. The forward and backward jump operators $\sigma,
\rho:\mathbb{T}\rightarrow\mathbb{T}$ and the graininess
$\mu:\mathbb{T}\rightarrow\mathbb{R}^+$ are defined, respectively,
by
\[
\sigma(t)=\inf\{s\in\mathbb{T}:s>t\},\,\,\,\,
\rho(t)=\sup\{s\in\mathbb{T}:s<t\}\,\,\,\,{\rm
and}\,\,\,\,\mu(t)=\sigma(t)-t.
\]

A point $t\in\mathbb{T}$ is called left-dense if $t>\inf\mathbb{T}$
and $\rho(t)=t$, left-scattered if $\rho(t)<t$, right-dense if
$t<\sup\mathbb{T}$ and $\sigma(t)=t$, and right-scattered if
$\sigma(t)>t$. If $\mathbb{T}$ has a left-scattered maximum $m$,
then $\mathbb{T}^k=\mathbb{T}\setminus\{m\}$; otherwise
$\mathbb{T}^k=\mathbb{T}$. If $\mathbb{T}$ has a right-scattered
minimum $m$, then $\mathbb{T}_k=\mathbb{T}\setminus\{m\}$; otherwise
$\mathbb{T}_k=\mathbb{T}$.
The time scale interval  $[a, b]_\mathbb{T}$  will
be defined by
$[a, b]_\mathbb{T} = \{t:  t\in [a,b]\cap \mathbb{T}\}$.

A function $f:\mathbb{T}\rightarrow\mathbb{R}$ is right-dense
continuous provided it is continuous at right-dense point in
$\mathbb{T}$ and its left-side limits exist at left-dense points in
$\mathbb{T}$. If $f$ is continuous at each right-dense point and
each left-dense point, then $f$ is said to be continuous function on
$\mathbb{T}$.

For $y:\mathbb{T}\rightarrow\mathbb{R}$ and $t\in\mathbb{T}^k$, we
define the delta derivative of $y(t)$, $y^\Delta(t)$, to be the
number (if it exists) with the property that for a given
$\varepsilon>0$, there exists a neighborhood $U$ of $t$ such that
\[
|[y(\sigma(t))-y(s)]-y^\Delta(t)[\sigma(t)-s]|<\varepsilon|\sigma(t)-s|
\]
for all $s\in U$.

If $y$ is continuous, then $y$ is right-dense continuous, and if $y$
is delta differentiable at $t$, then $y$ is continuous at $t$.

Let $y$ be right-dense continuous. If $Y^{\Delta}(t)=y(t)$, then we
define the delta integral by
$\int_a^{t}y(s)\Delta s=Y(t)-Y(a)$ if $a,t\in \mathbb{T}$; $\int_a^{t}y(s)\Delta s=Y(\rho(t))-Y(a)$ if $a\in \mathbb{T},t\notin \mathbb{T}, \rho(t)\geq a$ and
$\int_a^{t}y(s)\Delta s=Y(\sigma(t))-Y(a)$ if $a\in \mathbb{T},t\notin \mathbb{T}, \sigma(t)\leq a$; and so on.

\begin{definition}\cite{fink}\label{fk} A subset $S$ of $\mathbb{R}$ is called relatively dense or relatively dense in $\mathbb{R}$ if there exists a positive
number $L$ such that $[a, a + L] \cap S  \neq \emptyset$ for all $a\in \mathbb{R}$; a subset $S\subset A\subset \mathbb{R}$  is called relatively dense in $A$ if there exists a positive
number $L$ such that $[a, a + L]{\cap A} \cap S \neq \emptyset$ for all $a\in A$.The number $L$ is called the inclusion
length.
\end{definition}

Let $\mathbb{T}$ be a time scale and $\tau\in \mathbb{R}$, we denote $\mathbb{T}^\tau=\{t+\tau: t\in \mathbb{T}\}$.  
\begin{definition}\label{def21}\cite{li1,lm, k1}  A time scale $\mathbb{T}$ is called an invariant under a translation time scale  or a periodic time scale if
\begin{eqnarray}\label{amo}
\Pi=\big\{\tau\in\mathbb{R}: \mathbb{T}\cap \mathbb{T}^{\pm \tau}=\mathbb{T}\big\}\neq\{0\}.
\end{eqnarray}
\end{definition}

\begin{remark}\label{r21}
Since the set defined by \eqref{epe} is essentially the same set defined by \eqref{amo}, the almost periodic time scale under  Definition \ref{pe} is the periodic time scale under  Definition \ref{def21}.
If $\mathbb{T}$ is   invariant under a translation time scale, then $\sup \mathbb{T}=+\infty$ and $\inf \mathbb{T}=-\infty$.
\end{remark}

\begin{definition}\cite{lili1}\label{res1}
A time scale $\mathbb{T}$ is called an almost periodic time scale if
\[
\Pi:=\big\{\tau\in \mathbb{R}:\mathbb{T}_\tau\neq \emptyset\big\}
\]
is relatively dense in $\mathbb{R}$,  where $\mathbb{T}_\tau=\mathbb{T}\cap\{\mathbb{T}-\tau\}$ or $\mathbb{T}_\tau=\mathbb{T}\cap\{\mathbb{T}\pm\tau\}$.
\end{definition}

\begin{definition} \cite{lili3}\label{new2}
A time scale $\mathbb{T}$ is called an almost periodic time scale if the set $$\Pi:=\big\{\tau\in \mathbb{R}:\mathbb{T}_\tau\neq \emptyset\big\},$$
where $\mathbb{T}_\tau=\mathbb{T}\cap\{\mathbb{T}-\tau\}=\mathbb{T}\cap \{t-\tau: t\in \mathbb{T}\}$,  satisfies
\begin{itemize}
  \item [$(i)$] $\Pi\neq\{0\}$,
   \item [$(ii)$] if $\tau_1,\tau_2\in \Pi,$ then $\tau_1\pm \tau_2\in \Pi$.
\end{itemize}
\end{definition}
\begin{remark}\label{ar1}
If $\mathbb{T}$ is an almost periodic time scale under Definition \ref{def21}, then $\mathbb{T}$ is also an almost periodic time scale under Definitions \ref{res1} and \ref{new2}. If $\mathbb{T}$ is an almost periodic time scale under Definition \ref{new2}, then $\mathbb{T}$ is also an almost periodic time scale under Definition \ref{res1}.
\end{remark}
As a slightly modified version of Definition \ref{new2}, the following definition of almost periodic time scales is given in \cite{llm}:
\begin{definition}\cite{llm}\label{d1}
A time scale $\mathbb{T}$ is called an almost periodic time scale if the set $$\Pi:=\big\{\tau\in \mathbb{R}:\mathbb{T}_\tau\neq \emptyset\big\},$$
where $\mathbb{T}_\tau=\mathbb{T}\cap\{\mathbb{T}-\tau\}=\mathbb{T}\cap \{t-\tau: t\in \mathbb{T}\}$,  satisfies
\begin{itemize}
  \item [$(i)$] $\Pi\neq\{0\}$,
   \item [$(ii)$] if $\tau_1,\tau_2\in \Pi,$ then $\tau_1\pm \tau_2\in \Pi$,
  \item [$(iii)$]$\widetilde{\mathbb{T}}:=\mathbb{T}(\Pi)=\bigcap\limits_{\tau\in\Pi}\mathbb{T}_\tau\neq \emptyset$.
\end{itemize}
\end{definition}

\begin{definition}\label{defli2}\cite{li1}
Let $\mathbb{T}$ be an invariant under a translation
time scale.  A function  $f\in
C(\mathbb{T}\times D,\mathbb{E}^n)$ is  called an almost
periodic function in $t\in \mathbb{T}$ uniformly for $x\in D$ if the
$\varepsilon$-translation set of $f$
$$E\{\varepsilon,f,S\}=\{\tau\in\Pi:|f(t+\tau,x)-f(t,x)|<\varepsilon,\,\,
\forall (t,x)\in   \mathbb{T}\times S\}$$ is a relatively
dense set in $\mathbb{T}$ for all $\varepsilon>0$ and   for each
compact subset $S$ of $D$; that is, for any given $\varepsilon>0$
and each compact subset $S$ of $D$, there exists a constant
$l(\varepsilon,S)>0$ such that each interval of length
$l(\varepsilon,S)$ contains   a $\tau(\varepsilon,S)\in
E\{\varepsilon,f,S\}$ such that
\begin{equation*}
|f(t+\tau,x)-f(t,x)|<\varepsilon, \,\,\forall t\in
\mathbb{T}\times S.
\end{equation*}
This   $\tau$ is called the $\varepsilon$-translation number of $f$.
\end{definition}

\begin{lemma}\cite{mo}\label{mo1}
Let $\mathbb{T}$ be  an invariant under a translation time
scale. Then one has
\begin{itemize}
  \item [$(i)$] $\Pi\subset \mathbb{T}\Leftrightarrow 0\in \mathbb{T}$;
  \item [$(ii)$] $\Pi\cap \mathbb{T}=\emptyset \Leftrightarrow 0\notin \mathbb{T}$.
\end{itemize}
\end{lemma}

\section{Some remarks}
  \setcounter{equation}{0}
\indent

By Definition \ref{defli2}, we see that $E\{\varepsilon,f,S\}\subset \Pi$. So, according to Lemma \ref{mo1} and also, as pointed out in \cite{wa1}, $\mathbb{T}\cap\Pi=\emptyset$ if $0\notin \mathbb{T}$, and so, in this case,   $E\{\varepsilon,f,S\}$ can not be a relatively dense set in $\mathbb{T}$. To fix this flaw, the authors of \cite{wa1} use ``$E\{\varepsilon,f,S\}$  is a relatively dense set  in $\Pi$"  to replace ``$E\{\varepsilon,f,S\}$  is a relatively dense set in $\mathbb{T}$" in  Definition \ref{defli2} to give a  correction  of Definition \ref{defli2}, that is, Definition 1.8 in \cite{wa1}. Here, we will use ``$E\{\varepsilon,f,S\}$  is a relatively dense set in $\mathbb{R}$"  to replace ``$E\{\varepsilon,f,S\}$  is a relatively dense set in $\mathbb{T}$" in  Definition \ref{defli2} to give another correction of Definition \ref{defli2} as follows.

\begin{definition}\label{defli3}
Let $\mathbb{T}$ be an invariant under a translation
time scale. A function  $f\in
C(\mathbb{T}\times D,\mathbb{E}^n)$ is  called an almost
periodic function in $t\in \mathbb{T}$ uniformly for $x\in D$ if the
$\varepsilon$-translation set of $f$
$$E\{\varepsilon,f,S\}=\{\tau\in\Pi:|f(t+\tau,x)-f(t,x)|<\varepsilon,\,\,
\forall (t,x)\in   \mathbb{T}\times S\}$$ is a relatively
dense set in $\mathbb{R}$ for all $\varepsilon>0$ and   for each
compact subset $S$ of $D$; that is, for any given $\varepsilon>0$
and each compact subset $S$ of $D$, there exists a constant
$l(\varepsilon,S)>0$ such that each interval of length
$l(\varepsilon,S)$ contains   a $\tau(\varepsilon,S)\in
E\{\varepsilon,f,S\}$ such that
\begin{equation*}
|f(t+\tau,x)-f(t,x)|<\varepsilon, \,\,\forall t\in
\mathbb{T}\times S.
\end{equation*}
   $\tau$ is called the $\varepsilon$-translation number of $f$.
\end{definition}

The authors of \cite{wa1} also pointed out that there exist many integrals in \cite{li1} such as the following:
\[
\int_{t}^{t+l}f(s) \Delta s,
\]
but $t\in \mathbb{T}$ cannot guarantee that $t + l \in \mathbb{T}$, where $l \in \mathbb{R}$ is an inclusion length in Definition \ref{fk}. However, after we have thoroughly and carefully checked the whole paper \cite{li1},
there is no such integral in \cite{li1} at all.
Besides, even if there is such integral in paper \cite{li1}, according to our convention
\[
\int_{t}^{t+l}f(s) \Delta s:=\int_{t}^{\rho(t+l)}f(s) \Delta s,\quad \mathrm{if}\quad t\in \mathbb{T},\,\, t+l\notin \mathbb{T}.
\]
Hence, the integral $\int_{t}^{t+l}f(s) \Delta s$ is well defined.

\begin{remark}
From the above, we see that if we adopt   Definition \ref{defli3} as the definition of almost periodic functions on time scales, all the results of \cite{li1} remain true.
\end{remark}

In order to extend the concept of almost periodic  time scales, authors of \cite{wa1} give the following result:

\begin{theorem}(Theorem 2.4 in \cite{wa1})\label{wat1} Let $\mathbb{T}$ be an arbitrary time scale with $\sup \mathbb{T} = +\infty, \inf \mathbb{T} = -\infty$. If $\mu : \mathbb{T}\rightarrow\mathbb{R}^+$
is bounded, then $\mathbb{T}$ contains at least one invariant under the translation unit (that is, $\mathbb{T}$ contains at least one   sub time scale that is  an invariant under a translation time scale).
\end{theorem}

Unfortunately, the following example shows that Theorem \ref{wat1} is incorrect.

\begin{example}\label{ex1}
Take $\mathbb{T}=\{-2k,2k+1:k\in \mathbb{N}\}$, then $\sup\limits_{t\in \mathbb{T}} \mu(t)=5$. According to Remark \ref{r21}, if $\mathbb{T}_*$ is a sub time scale of $\mathbb{T}$ that is  an invariant under a translation time scale, then $\sup \mathbb{T}_*=+\infty$ and $\inf \mathbb{T}_*=-\infty$.

Let $\mathbb{\hat{T}}$ be   any subset of $\mathbb{T}$ with $\sup \mathbb{\hat{T}}=+\infty$ and $\inf \mathbb{\hat{T}}=-\infty$, then $\mathbb{\hat{T}}$  can be expressed as
\[
\mathbb{\hat{T}}=\{-2k_{n},2l_m+1: n,m\in \mathbb{N}\},
\]
  where $\{k_n\}_{n=1}^{\infty},\{l_m\}_{m=1}^{\infty}\subset \mathbb{N}, k_n< k_{n+1}, l_m< l_{m+1}$ and $\lim\limits_{n\rightarrow\infty}k_n=\lim\limits_{m\rightarrow\infty}l_m=+\infty$. Then for every $\tau\in \mathbb{R}$,
$$\mathbb{\hat{T}}^\tau=\{-2k_{n}+\tau,2l_m+1+\tau: n,m\in \mathbb{N}\}.$$

In the following, we divide  five cases to show that for every $\tau \in \mathbb{R}\setminus\{0\}$, $\mathbb{\hat{T}}\cap \mathbb{\hat{T}}^\tau\neq \mathbb{\hat{T}}$, that is, $\mathbb{T}$ contains no any sub time scale that is  an invariant under a translation time scale.

 Case 1: It is easy to see that $\mathbb{T}\cap \mathbb{\hat{T}}^\tau=\emptyset$ for $\tau\in \mathbb{R}\setminus\mathbb{Z}$.

   Case 2: If $\tau$ is a positive even number, then  $2l_{1}+1\in \mathbb{\hat{T}}$ and $2l_{1}+1\notin \mathbb{\hat{T}}^\tau$, hence   $\mathbb{\hat{T}}\cap \mathbb{\hat{T}}^\tau\neq \mathbb{\hat{T}}$.

  Case 3: If $\tau$ is a negative even number, then  $-2k_1\in \mathbb{\hat{T}}$ and $-2k_1\notin \mathbb{\hat{T}}^\tau$, hence   $\mathbb{\hat{T}}\cap \mathbb{\hat{T}}^\tau\neq \mathbb{\hat{T}}$.

  Case 4: If $\tau$ is a positive odd number, then  there are at last only finite many positive odd numbers in $\{-2k_n+\tau\}_{n=1}^{\infty}$, that is, $\mathbb{\hat{T}}^\tau$ only contains at last finite many positive odd numbers. Since there are infinite many positive odd numbers in $\mathbb{\hat{T}}$, hence   $\mathbb{\hat{T}}\cap \mathbb{\hat{T}}^\tau\neq \mathbb{\hat{T}}$.

  Case 5: If $\tau$ is a negative odd number, then  there are at last only finite many negative even numbers in $\{2l_m+1+\tau\}_{n=1}^{\infty}$, that is, $\mathbb{\hat{T}}^\tau$ only contains at last finite many negative even numbers. Since there are infinite many negative even numbers in $\mathbb{\hat{T}}$, hence   $\mathbb{\hat{T}}\cap \mathbb{\hat{T}}^\tau\neq \mathbb{\hat{T}}$.
\end{example}
\begin{remark}\label{r32}
Example \ref{ex1} shows that there exists a time scale  that satisfies all the conditions of Theorem \ref{wat1}, but it    contains no sub time scale that is an invariant under a translation time scale.
Therefore,  Theorem \ref{wat1} is incorrect.
\end{remark}
\begin{remark}
The set $\widetilde{\mathbb{T}}$ of Definition \ref{d1} is the sub-invariant under
the translation unit in $\mathbb{T}$ of Definitions 2.2 and 2.3 in \cite{wa1}.
\end{remark}

In \cite{wa2}, the following definitions are given:
\begin{definition}(Definition 2.4 in \cite{wa2})\label{e24}
Let $\mathbb{T}$ be a time scale, we say $\mathbb{T}$ is a zero-periodic time scale if and only if
there exists no nonzero real number $\omega$ such that $t+\omega\in\mathbb{T}$ for all $t\in\mathbb{T}$.
\end{definition}

\begin{definition}(Definition 2.5 in \cite{wa2})\label{e25}
A time scale sequence $\{\mathbb{T}_i\}_{i\in \mathbb{Z}^+}$ is called well-connected if and only if for $i\neq j$, we have $\mathbb{T}_i\cap\mathbb{T}_j=\{t_{ij}^{k}\}_{k\in \mathbb{Z}}$, where $\{t_{ij}^{k}\}$ is a countable points set or an empty set, and $t_{ij}^{k}$ is called the connected point between $\mathbb{T}_i$ and $\mathbb{T}_j$ for each $k\in \mathbb{Z}$, the set $\{t_{ij}^{k}\}$ is called the
connected points set of this well-connected sequence.
\end{definition}

\begin{definition}(Definition 2.6 in \cite{wa2})\label{e26}
Let $\mathbb{T}$ be an infinite time scale. We say $\mathbb{T}$ is a changing-periodic or a
piecewise-periodic time scale if the following conditions are fulfilled:
\begin{enumerate}
  \item [$(a)$]
$\mathbb{T}=(\bigcup_{i=1}^{\infty}\mathbb{T}_i)\cup\mathbb{T}_r$ and $\{\mathbb{T}_i\}_{i\in \mathbb{Z}^+}$ is a well-connected time scale sequence, where
$\mathbb{T}_r=\bigcup_{i=1}^{k}[\alpha_i,\beta_i]$ and $k$ is some finite number, and $[\alpha_i,\beta_i]$ are closed intervals for $i=1,2,\ldots,k$
or $\mathbb{T}_r=\emptyset$;
\item [$(b)$]
$S_i$ is a nonempty subset of $\mathbb{R}$ with $0\notin S_i$ for each $i\in \mathbb{Z}^+$ and $\Pi=(\bigcup_{i=1}^{\infty}S_i)\cup R_0$, where $R_0=\{0\}$ or $R_0=\emptyset$;
\item [$(c)$]
for all $t\in\mathbb{T}_i$ and all $\omega\in S_i$, we have $t+\omega\in\mathbb{T}_i$, i.e., $\mathbb{T}_i$ is an $\omega$-periodic time scale;
\item [$(d)$]
for $i\neq j$, for all $t\in\mathbb{T}_i\setminus\{t_{ij}^{k}\}$ and all $\omega\in S_j$, we have $t+\omega\notin\mathbb{T}$, where $\{t_{ij}^{k}\}$ is the
connected points set of the timescale sequence $\{\mathbb{T}_i\}_{i\in \mathbb{Z}^+}$;
\item [$(e)$]
$R_0=\{0\}$ if and only if $\mathbb{T}_r$ is a zero-periodic time scale and $R_0=\emptyset$ if and only if $\mathbb{T}_r=\emptyset$;
\end{enumerate}
and the set $\Pi$ is called a changing-periods set of $\mathbb{T}$, $\mathbb{T}_i$ is called the periodic sub-timescale
of $\mathbb{T}$ and $S_i$ is called the periods subset of $\mathbb{T}$ or the periods set of $\mathbb{T}_i$, $\mathbb{T}_r$ is called the remain
timescale of $\mathbb{T}$ and $R_0$ the remain periods set of $\mathbb{T}$.
\end{definition}

In \cite{wa2}, the following results are given:
\begin{theorem}(Theorem 2.11 in \cite{wa2}) \label{wat2} If $\mathbb{T}$ is an infinite time scale and the graininess function $\mu : t \rightarrow \mathbb{R}^+$ is
bounded, then $\mathbb{T}$ is a changing-periodic time scale.
\end{theorem}

\begin{theorem}(Theorem 2.21 in \cite{wa2})\label{wat3}
(Decomposition theorem of time scales) Let $\mathbb{T}$ be an infinite time scale
and the graininess function $\mu : \mathbb{T}\rightarrow \mathbb{R}^+$ be bounded, then $\mathbb{T}$ is a changing-periodic time
scale, i.e., there exists a countable periodic decomposition such that $\mathbb{T}=(\bigcup_{i=1}^{\infty}\mathbb{T}_i)\cup\mathbb{T}_r$ and $\mathbb{T}_i$ is an $\omega$-periodic sub-timescale, $\omega\in S_i, i \in \mathbb{Z}^+$, where $\mathbb{T}_i, S_i, \mathbb{T}_r$ satisfy the conditions in
Definition \ref{e26}.
\end{theorem}

\begin{theorem}(Theorem 2.23 in \cite{wa2}) \label{wat4}
(Periodic coverage theorem of time scales) Let $\mathbb{T}$ be an infinite time scale
and the graininess function $\mu: \mathbb{T}\rightarrow\mathbb{R}^+$
be bounded, then $\mathbb{T}$ can be covered by countable periodic time scales.
\end{theorem}

Unfortunately, we have the following remark:
\begin{remark}
Theorems \ref{wat2} and \ref{wat3} are incorrect. The reasons are as follows:

According to Definition \ref{e26}, if $\mathbb{T}$ is a changing-periodic time scale, then $\mathbb{T}$ can be expressed as $\mathbb{T}=(\bigcup_{i=1}^{\infty}\mathbb{T}_i)\cup\mathbb{T}_r$, where $\mathbb{T}_i\subset \mathbb{T}$ is a periodic time scale. Hence, Theorem \ref{wat2} is essentially equivalent to Theorem \ref{wat1}. Therefore,  in view of Remark \ref{r32}, Theorem \ref{wat2} is incorrect. For the same reason, Theorem \ref{wat3} is also incorrect.

As to Theorem \ref{wat4}, since its proof is based on Theorem \ref{wat2}, we cannot judge  the correctness of it.
\end{remark}

The following concept of an index function for changing-periodic time scales in \cite{wa2} has flaws.

\begin{definition}(Definition 2.8 in \cite{wa2})\label{e28}
Let $\mathbb{T}$ be a changing-periodic time scale, then the function $\tau$
\begin{eqnarray*}
\tau&:&\mathbb{T}\mapsto\mathbb{Z}^+\cup\{0\},\\
&&\bigg(\bigcup_{i=1}^{\infty}\mathbb{T}_i\bigg)\mapsto i,\,\, where\,\, t\in\mathbb{T}_i,\,\,i\in \mathbb{Z}^+,\\
&&\mathbb{T}_r\mapsto 0,\,\, where\,\, t\in\mathbb{T}_r,\\
&&t\mapsto\tau_t
\end{eqnarray*}
is called an index function for $\mathbb{T}$, where the corresponding periods set of $\mathbb{T}_{\tau_t}$ is denoted as $S_{\tau_t}$. In what follows we shall call $S_{\tau_t}$ the adaption set generated by $t$, and all the elements in $S_{\tau_t}$ will be called the adaption factors for $t$.
\end{definition}

\begin{remark}
According to  Definition \ref{e25} and Definition \ref{e26}, for $i\neq j$, we have $\mathbb{T}_i\cap\mathbb{T}_j=\{t_{ij}^{k}\}_{k\in \mathbb{Z}}$, where $\{t_{ij}^{k}\}$ is a countable points set or an empty set, and
$\mathbb{T}_r\cap \mathbb{T}_i$ may not be an empty set. Therefore, the concept of an index function for changing-periodic time scales is not well defined.
\end{remark}
\begin{definition}(Definition 3.1 in \cite{wa2})\label{e31}
Let $\mathbb{T}$ be a changing-periodic time scale, i.e., $\mathbb{T}$ satisfies Definition \ref{e26}. A function $f\in C(\mathbb{T}\times D, E^n)$ is called a local-almost periodic function in $t\in\mathbb{T}$ uniformly for $x\in D$ if the $\varepsilon$-translation numbers set of $f$,
\begin{eqnarray*}
E\{\varepsilon, f, S\}=\{\tilde{\tau}\in S_{\tau_t}:|f(t+\tilde{\tau},x)-f(t,x)|<\varepsilon\,\, for\,\, all \,\,(t,x)\in\mathbb{T}\times S\}
\end{eqnarray*}
is a relatively dense set for all $\varepsilon>0$ and for each compact subset $S$ of $D$; that is, for any given $\varepsilon>0$ and each compact subset $S$ of $D$, there exists a constant $l(\varepsilon, S)>0$ such that each interval of length $l(\varepsilon,S)$ contains $\tilde{\tau}(\varepsilon, S)\in E\{\varepsilon, f, S\}$ such that
\begin{eqnarray*}
|f(t+\tilde{\tau},x)-f(t,x)|<\varepsilon\,\, for\,\, all \,\,(t,x)\in\mathbb{T}\times S;
\end{eqnarray*}
here, $\tilde{\tau}$ is called the $\varepsilon$-local translation number of $f$ and $l(\varepsilon,S)$ is called the local inclusion
length of $E\{\varepsilon, f, S\}$.
\end{definition}

\begin{definition}(Definition 3.2 in \cite{wa2})\label{e32}
Assume that $\mathbb{T}$ is a changing-periodic time scale. Let $f(t,x)\in C(\mathbb{T}\times D, E^n)$ if for any given adaption factors sequence $(\alpha^\tau)'\subset S_{\tau_t}$, there exists a subsequence $\alpha^\tau\subset(\alpha^\tau)'$ such that $T_{\alpha^\tau}f(t,x)$ exists uniformly on $\mathbb{T}\times S$, then $f(t,x)$ is called a local-almost periodic function in $t$ uniformly for $x\in D$.
\end{definition}

\begin{remark}
Since Definition \ref{e31} and Definition \ref{e32} are based on the concept of an index function for changing-periodic time scales, they are not well defined.
\end{remark}

In \cite{wa3,wa4}, the following definitions of almost periodic time scales and almost periodic functions on time scales are given:
\begin{definition}\cite{wa4, wa3}\label{e52}
We say that $\mathbb{T}$ is an almost periodic time scale if for any given $\varepsilon_1>0$, there exists a constant $l(\varepsilon_1)>0$ such that each interval of length $l(\varepsilon_1)$ contains a $\tau(\varepsilon_1)$ such that
\begin{eqnarray*}
d(\mathbb{T},\mathbb{T}^\tau)<\varepsilon_1;
\end{eqnarray*}
that is, for any $\varepsilon_1>0$, the following set
\begin{eqnarray*}
E\{\mathbb{T},\varepsilon_1\}=\{\tau\in\mathbb{R}:d(\mathbb{T}^\tau,\mathbb{T})<\varepsilon_1\}
\end{eqnarray*}
is relatively dense. This $\tau$ is called the $\varepsilon_1$-translation number of $\mathbb{T}$, $l(\varepsilon_1)$ is called the inclusion length of $E\{\mathbb{T},\varepsilon_1\}$, and $E\{\mathbb{T},\varepsilon_1\}$ is called the $\varepsilon_1$-translation set of $\mathbb{T}$.
\end{definition}

In the following, we  denote
$\Pi_\varepsilon=E\{\mathbb{T},\varepsilon\}$ and $\mathbb{T}^{\Pi_\varepsilon}=\{\mathbb{T}^{-\tau}:-\tau\in E\{\mathbb{T},\varepsilon\}\}.$

\begin{definition}(Definition 5.5 in \cite{wa4})\label{e55}
Let $\mathbb{T}$ be an almost periodic time scale under Definition \ref{e52}. A function $f\in C(\mathbb{T}\times D, E^n)$ is called an almost periodic function in $t\in\mathbb{T}$ uniformly for $x\in D$ if the $\varepsilon_2$-translation set of $f$
\begin{eqnarray*}
E\{\varepsilon_2, f, S\}=\{\tau\in\Pi_{\varepsilon_1}:|f(t+\tau,x)-f(t,x)|<\varepsilon_2, 
 \,\, for\,\, all \,\,(t,x)\in(\mathbb{T}\cap(\cup_{-\tau}\mathbb{T}^{\Pi_{\varepsilon_{1}}}))\times S\}
\end{eqnarray*}
is a relatively dense set for all $\varepsilon_2>\varepsilon_1>0$ and for each compact subset $S$ of $D$; that is, for any given $\varepsilon_2>\varepsilon_1>0$ and each compact subset $S$ of $D$, there exists a constant $l(\varepsilon_2, S)>0$ such that each interval of length $l(\varepsilon_2,S)$ contains $\tau(\varepsilon_2, S)\in E\{\varepsilon_2, f, S\}$ such that
\begin{eqnarray*}
|f(t+\tau,x)-f(t,x)|<\varepsilon_2,\,\, for\,\, all \,\,(t,x)\in(\mathbb{T}\cap(\cup_{-\tau}\mathbb{T}^{\Pi_{\varepsilon_{1}}}))\times S.
\end{eqnarray*}
This $\tau$ is called the $\varepsilon_2$-translation number of $f$ and $l(\varepsilon_2,S)$ is called the inclusion
length of $E\{\varepsilon_2, f, S\}$.
\end{definition}
\begin{remark}Since the fact that $\mathbb{T}$ is an almost periodic time scale under Definition \ref{e52} may do  not guarantee that the set $\{\tau \in \Pi_{\varepsilon_1}=E\{\mathbb{T},\varepsilon_1\}: \mathbb{T}\cap \mathbb{T}^\tau\neq \emptyset\}$ is relatively dense. Therefore,  Definition \ref{e55} is not well defined. An correction for this, we refer to   \cite{lili1}.
\end{remark}
\begin{definition}(Definition 6.2 in \cite{wa4})\label{e62}
Let $\mathbb{T}$ be an almost periodic time scale. A function $f\in C(\mathbb{T}\times D, E^n)$ is called an $\varepsilon^*$-local almost periodic function in $t\in\mathbb{T}$ uniformly for $x\in D$  if there exists some fixed
$\varepsilon^*>0$, such that for all $\varepsilon_1>\varepsilon^*$ and $\varepsilon^*<d(\mathbb{T},\mathbb{T}^\tau)<\varepsilon_1$,
the $\varepsilon_2$-translation set of $f$
\begin{eqnarray*}
E\{\varepsilon_2, f, S\}=\{\tau\in\Pi_{\varepsilon_1}:|f(t+\tau,x)-f(t,x)|<\varepsilon_2,
\,\, for\,\, all \,\,(t,x)\in(\mathbb{T}\cap \mathbb{T}^{-\tau})\times S\}
\end{eqnarray*}
is a relatively dense set for all $\varepsilon_2>\varepsilon_1>0$ and for each compact subset $S$ of $D$; that is, for any given $\varepsilon_2>0$ and each compact subset $S$ of $D$, there exists a constant $l(\varepsilon_2, S)>0$ such that each interval of length $l(\varepsilon_2,S)$ contains a $\tau(\varepsilon_2, S)\in E\{\varepsilon_2, f, S\}$ such that
\begin{eqnarray*}
|f(t+\tau,x)-f(t,x)|<\varepsilon_2,\,\, for\,\, all \,\,(t,x)\in(\mathbb{T}\cap\mathbb{T}^{-\tau})\times S.
\end{eqnarray*}
This $\tau$ is called the $\varepsilon_2$-translation number of $f$ and $l(\varepsilon_2,S)$ is called the inclusion
length of $E\{\varepsilon_2, f, S\}$.
\end{definition}

\begin{remark}If we take $\mathbb{T}=\mathbb{Z}$, then $\mathbb{T}\cap \mathbb{T}^\tau=\emptyset$ for $\tau\in \{\tau\in \mathbb{R}: \varepsilon^*<d(\mathbb{T},\mathbb{T}^\tau)<\varepsilon_1<1\}$.
Therefore, Definition \ref{e62} is not well defined.
\end{remark}

As a closing of this note, we shall point out that, in order to define  local-almost periodic functions, local-almost automorphic functions and so on,  on time scales, one needs a proper definition of almost periodic time scales, which can support these classes of functions on time scales,  the almost periodic time scale under Definition \ref{res1} may be the most general almost periodic time scale. However, since a time scale under Definition \ref{new2} has the concrete property $(ii)$ of Definition \ref{new2},  for conveniens, the almost periodic time scale under Definition \ref{new2} may be the best one on which one can define and investigate local-almost periodic functions, local-almost automorphic functions and so on.\\

\textbf{Conflict of Interests}\\

The author  declares that there is no conflict of interests
regarding the publication of this paper.

\end{document}